\documentclass{article}
\usepackage{amssymb}
\usepackage{amsmath}

\newtheorem{theorem}{Theorem}

\newtheorem{definition}[theorem]{Definition}

\newtheorem{proposition}[theorem]{Proposition}
\newtheorem{remark}{Remark}

\begin{document}

\title{\textbf{Pointwise strong approximation of almost periodic
functions }}
\author{\textbf{Rados\l awa Kranz, W\l odzimierz \L enski \ and Bogdan Szal} 
\\
University of Zielona G\'{o}ra\\
Faculty of Mathematics, Computer Science and Econometrics\\
65-516 Zielona G\'{o}ra, ul. Szafrana 4a\\
P O L A N D\\
R.Kranz@wmie.uz.zgora.pl, W.Lenski@wmie.uz.zgora.pl, \\
B.Szal @wmie.uz.zgora.pl}
\date{}
\maketitle

\begin{abstract}
We consider the class $GM\left( _{2}\beta \right) $ in pointwise estimate of
the deviations in strong mean of almost periodic functions from matrix means
of partial sums of their Fourier series.

\ \ \ \ \ \ \ \ \ \ \ \ \ \ \ \ \ \ \ \ 


\ \ \ \ \ \ \ \ \ \ \ \ \ \ \ \ \ \ \ 

\end{abstract}

\section{Introduction}

Let $S^{p}\;\left( 1<p\leq \infty \right) $ be the class of all almost
periodic functions in the sense of Stepanov with the norm 
\begin{equation*}
\Vert f\Vert _{S^{p}}:=\left\{ 
\begin{array}{c}
\sup\limits_{u}\left\{ \frac{1}{\pi }\int_{u}^{u+\pi }\mid f(t)\mid
^{p}dt\right\} ^{1/p}\text{ \ \ when \ \ }1<p<\infty \\ 
\sup\limits_{u}\mid f(u)\mid \text{ \ \ when \ \ }p=\infty .%
\end{array}%
\right.
\end{equation*}%
Suppose that the Fourier series of $f\in S^{p}$ has the form 
\begin{equation*}
Sf\left( x\right) =\sum_{\nu =-\infty }^{\infty }A_{\nu }\left( f\right)
e^{i\lambda _{\nu }x},\text{ \ \ where \ }A_{\nu }\left( f\right)
=\lim_{L\rightarrow \infty }\frac{1}{L}\int_{0}^{L}f(t)e^{-i\lambda _{\nu
}t}dt,
\end{equation*}%
with the partial sums\ 
\begin{equation*}
S_{\gamma _{k}}f\left( x\right) =\sum_{\left\vert \lambda _{\nu }\right\vert
\leq \gamma _{k}}A_{\nu }\left( f\right) e^{i\lambda _{\nu }x}
\end{equation*}%
and that $0=\lambda _{0}<\lambda _{\nu }<\lambda _{\nu +1}$ if $\nu \in 
\mathbb{N}
=\left\{ 1,2,3...\right\} ,$ $\underset{v\rightarrow \infty }{\lim }\lambda
_{\nu }=\infty ,$ $\lambda _{-\nu }=-\lambda _{\nu ,}$\ $\left\vert A_{\nu
}\right\vert +\left\vert A_{-\nu }\right\vert >0.$ Let\ $\Omega _{\alpha ,p}$
, with some fixed positive $\alpha $ , be the set of functions of class $%
S^{p}$ bounded on $U=\left( -\infty ,\infty \right) $ whose Fourier
exponents satisfy the condition 
\begin{equation*}
\lambda _{\nu +1}-\lambda _{\nu }\geq \alpha \text{ \ \ }\left( \nu \in 
\mathbb{N}
\right) .
\end{equation*}%
In case $f\in \Omega _{\alpha ,p}$\ 
\begin{equation*}
S_{\lambda _{k}}f\left( x\right) =\int_{0}^{\infty }\left\{ f\left(
x+t\right) +f\left( x-t\right) \right\} \Psi _{\lambda _{k},\lambda
_{k}+\alpha }\left( t\right) dt,
\end{equation*}%
where 
\begin{equation*}
\Psi _{\lambda ,\eta }\left( t\right) =\frac{2\sin \frac{\left( \eta
-\lambda \right) t}{2}\sin \frac{\left( \eta +\lambda \right) t}{2}}{\pi
\left( \eta -\lambda \right) t^{2}}\text{ \ \ }\left( 0<\lambda <\eta ,\text{
\ }\left\vert t\right\vert >0\right) .
\end{equation*}

Let $A:=\left( a_{n,k}\right) $ be an infinite matrix of real nonnegative
numbers such that 
\begin{equation}
\sum_{k=0}^{\infty }a_{n,k}=1\text{, where }n=0,1,2,...\text{\ .}  \label{T1}
\end{equation}

Let us consider the strong mean 
\begin{equation}
H_{n,A,\gamma }^{q}f\left( x\right) =\left\{ \sum_{k=0}^{\infty
}a_{n,k}\left\vert S_{\gamma _{k}}f\left( x\right) -f\left( x\right)
\right\vert ^{q}\right\} ^{1/q}\text{ \ \ \ }\left( q>0\right) \text{.}
\label{S1}
\end{equation}%
As measures of approximation by the quantity (\ref{S1}), we use the best
approximation of $f$ by entire functions $g_{\sigma }$ of exponential type $%
\sigma $\ bounded on the real axis, shortly $g_{\sigma }\in B_{\sigma }$ and
the moduli of continuity of $\ f$ defined by the formulas%
\begin{equation*}
E_{\sigma }(f)_{S^{p}}=\inf_{g_{\sigma }}\left\Vert f-g_{\sigma }\right\Vert
_{S^{p}},
\end{equation*}%
\begin{equation*}
\omega f\left( \delta \right) _{X}=\sup_{\left\vert t\right\vert \leq \delta
}\left\Vert f\left( \cdot +t\right) -f\left( \cdot \right) \right\Vert _{X%
\text{ }},\text{ \ \ }X=C_{2\pi }\text{ or }X=S^{p}
\end{equation*}%
and%
\begin{equation*}
w_{x}f(\delta )_{S^{p}}:=\left\{ \frac{1}{\delta }\int_{0}^{\delta
}\left\vert \varphi _{x}\left( t\right) \right\vert ^{p}dt\right\} ^{1/p},
\end{equation*}%
respectively.

Recently, L. Leindler \cite{3} defined a new class of sequences named as
sequences of rest bounded variation, briefly denoted by $RBVS$, i.e. 
\begin{equation}
RBVS=\left\{ a:=\left( a_{n}\right) \in U:\sum\limits_{k=m}^{\infty
}\left\vert a_{k}-a_{k+1}\right\vert \leq K\left( a\right) \left\vert
a_{m}\right\vert \text{ for all }m\in U\right\} ,  \label{1}
\end{equation}%
where here and throughout the paper $K\left( a\right) $ always indicates a
constant depending only on $a$.

Denote by $MS$ the class of nonnegative and nonincreasing sequences

and by $GM$ the class of general monotone coefficients defined as follows (
see \cite{10}): 
\begin{equation}
GM=\left\{ a:=\left( a_{n}\right) \in U:\sum\limits_{k=m}^{2m-1}\left\vert
a_{k}-a_{k+1}\right\vert \leq K\left( a\right) \left\vert a_{m}\right\vert 
\text{ for all }m\in U\right\} .  \label{3}
\end{equation}%
Then it is obvious that 
\begin{equation*}
MS\subset RBVS\subset GM\text{.}
\end{equation*}%
In \cite{6, 10, 11, 12} was defined the class of $\beta -$general monotone
sequences as follows:

\begin{definition}
Let $\beta :=\left( \beta _{n}\right) $ be a nonnegative sequence. The
sequence of complex numbers $a:=\left( a_{n}\right) $ is said to be $\beta -$%
general monotone, or $a\in GM\left( \beta \right) $, if the relation 
\begin{equation}
\sum\limits_{k=m}^{2m-1}\left| a_{k}-a_{k+1}\right| \leq K\left( a\right)
\beta _{m}  \label{4}
\end{equation}
holds for all $m$.
\end{definition}

In the paper \cite{12} Tikhonov considered, among others, the following
examples of the sequences $\beta _{n}:$

(1) $_{1}\beta _{n}=\left\vert a_{n}\right\vert ,$

(2) $_{2}\beta _{n}=\sum\limits_{k=\left[ n/c\right] }^{\left[ cn\right] }%
\frac{\left\vert a_{k}\right\vert }{k}$ for some $c>1$.

It is clear that (see \cite[Remark 2.1]{12}) $GM\left( _{1}\beta \right) =GM$
and 
\begin{equation*}
GM\left( _{1}\beta +_{2}\beta \right) \equiv GM\left( _{2}\beta \right) .
\end{equation*}

Moreover, we assume that the sequence $\left( K\left( \alpha _{n}\right)
\right) _{n=0}^{\infty }$ is bounded, that is, that there exists a constant $%
K$ such that 
\begin{equation*}
0\leq K\left( \alpha _{n}\right) \leq K
\end{equation*}%
holds for all $n$, where $K\left( \alpha _{n}\right) $ denote the sequence
of constants appearing in the inequalities (\ref{1})-(\ref{4}) for the
sequences $\alpha _{n}:=\left( a_{n,k}\right) _{k=0}^{\infty }$.

Now we can give the conditions to be used later on. We assume that for all $n
$%
\begin{equation}
\sum\limits_{k=m}^{2m-1}\left\vert a_{n,k}-a_{n,k+1}\right\vert \leq
K\sum\limits_{k=[m/c]}^{[cm]}\frac{a_{n,k}}{k}  \label{5}
\end{equation}%
holds if $\alpha _{n}=\left( a_{n,k}\right) _{k=0}^{\infty }$ belongs to $%
GM\left( _{2}\beta \right) $, for $n=0,1,2,...$

We have shown in \cite{WLBS 2} the following theorem:

\begin{theorem}
If $f\in \Omega _{\alpha ,p},$ $p\geq q$, $\left( a_{n,k}\right)
_{k=0}^{\infty }\in GM\left( _{2}\beta \right) $ for all $n$, (\ref{T1}) and 
$\lim_{n\rightarrow \infty }a_{n,0}=0$ hold, then 
\begin{equation*}
\left\Vert H_{n,A,\gamma }^{q}f\right\Vert _{S^{p}}\ll \left\{
\sum_{k=0}^{\infty }a_{n,k}\omega ^{q}f\left( \frac{\pi }{k+1}\right)
_{S^{p}}\right\} ^{1/q},
\end{equation*}%
for $n=0,1,2,...$, where\ $\gamma =\left( \gamma _{k}\right) $ is a sequence
with $\gamma _{k}=\frac{\alpha k}{2}$.
\end{theorem}

In this paper we consider the class $GM\left( _{2}\beta \right) $ in
pointwise estimate of the quantity $H_{n,A,\gamma }^{q}f.$ Thus we present
some analog of the following result of P. Pych-Taberska\ (see \cite[Theorem 5%
]{PT}):

\begin{theorem}
If $f\in \Omega _{\alpha ,\infty }$ and $q\geq 2,$ then 
\begin{equation*}
\left\Vert H_{n,A,\gamma }^{q}f\right\Vert _{S^{\infty }}\ll \left\{ \frac{1%
}{n+1}\sum_{k=0}^{n}\left[ \omega f\left( \frac{\pi }{k+1}\right)
_{S^{\infty }}\right] ^{q}\right\} ^{1/q}+\frac{\left\Vert f\right\Vert
_{S^{\infty }}}{\left( n+1\right) ^{1/q}},
\end{equation*}%
for $n=0,1,2,...$, where $\gamma =\left( \gamma _{k}\right) $ is a sequence
with $\gamma _{k}=\frac{\alpha k}{2},$ $a_{n,k}=\frac{1}{n+1}$ when $k\leq n$
and $a_{n,k}=0$ otherwise.
\end{theorem}

We shall write $I_{1}\ll I_{2}$ if there exists a positive constant $K$ ,
sometimes depended on some parameters, such that $I_{1}\leq KI_{2}$.

\section{Statement of the results}

Let us consider a function $w_{x}$ of modulus of continuity type on the
interval $[0,+\infty ),$ i.e. a nondecreasing continuous function having the
following properties:\ $w_{x}\left( 0\right) =0,$ $w_{x}\left( \delta
_{1}+\delta _{2}\right) \leq w_{x}\left( \delta _{1}\right) +w_{x}\left(
\delta _{2}\right) $ for any $\delta _{1},\delta _{2}\geq 0$\ with $x$ such
that the set%
\begin{eqnarray*}
\Omega _{\alpha ,p}\left( w_{x}\right) &=&\left\{ f\in \Omega _{\alpha ,p}: 
\left[ \frac{1}{\delta }\int_{0}^{\delta }\left\vert \varphi _{x}\left(
t\right) -\varphi _{x}\left( t\pm \gamma \right) \right\vert ^{p}dt\right]
^{1/p}\ll w_{x}\left( \gamma \right) \right. \\
&&\left. \text{\ and \ }w_{x}f\left( \delta \right) _{p}\ll w_{x}\left(
\delta \right) \text{ \ , \ where \ }\gamma ,\delta >0\right\}
\end{eqnarray*}%
is nonempty. It is clear that $\Omega _{\alpha ,p}\left( w_{x}\right)
\subseteq \Omega _{\alpha ,p^{\prime }}\left( w_{x}\right) ,$ for $p^{\prime
}\leq p.$

We start with proposition

\begin{proposition}
If $f\in \Omega _{\alpha ,p}\left( w_{x}\right) $ and $q>0,$ then 
\begin{equation*}
\left\{ \frac{1}{n+1}\sum_{k=n}^{2n}\left\vert S_{\frac{\alpha k}{2}}f\left(
x\right) -f\left( x\right) \right\vert ^{q}\right\} ^{1/q}\ll w_{x}\left( 
\frac{\pi }{n+1}\right) +E_{\alpha n/2}\left( f\right) _{S^{p}},
\end{equation*}%
for $n=0,1,2,...$
\end{proposition}

Our main results are following

\begin{theorem}
If $f\in \Omega _{\alpha ,p}\left( w_{x}\right) ,$ $q>0$, $\left(
a_{n,k}\right) _{k=0}^{\infty }\in GM\left( _{2}\beta \right) $ for all $n$,
(\ref{T1}) and $\lim_{n\rightarrow \infty }a_{n,0}=0$ hold, then 
\begin{equation*}
H_{n,A,\gamma }^{q}f\left( x\right) \ll \left\{ \sum_{k=0}^{\infty }a_{n,k} 
\left[ w_{x}(\frac{\pi }{k+1})+E_{\frac{\alpha k}{2^{1+\left[ c\right] }}%
}\left( f\right) _{S^{p}}\right] ^{q}\right\} ^{1/q}
\end{equation*}%
for some $c>1$ and $n=0,1,2,...$, where\ $\gamma =\left( \gamma _{k}\right) $
is a sequence with $\gamma _{k}=\frac{\alpha k}{2}.$
\end{theorem}

\begin{theorem}
If $f\in \Omega _{\alpha ,p}\left( w_{x}\right) ,$ $q>0$, $\left(
a_{n,k}\right) _{k=0}^{\infty }\in MS$ for all $n$, (\ref{T1}) and $%
\lim_{n\rightarrow \infty }a_{n,0}=0$ hold, then 
\begin{equation*}
H_{n,A,\gamma }^{q}f\left( x\right) \ll \left\{ \sum_{k=0}^{\infty }a_{n,k} 
\left[ w_{x}(\frac{\pi }{k+1})+E_{\frac{\alpha k}{2}}\left( f\right) _{S^{p}}%
\right] ^{q}\right\} ^{1/q}
\end{equation*}%
for $n=0,1,2,...$, where\ $\gamma =\left( \gamma _{k}\right) $ is a sequence
with $\gamma _{k}=\frac{\alpha k}{2}$.
\end{theorem}

\begin{remark}
Since 
\begin{equation*}
\left\Vert \left[ \frac{1}{\delta }\int_{0}^{\delta }\left\vert \varphi
_{\cdot }\left( t\right) -\varphi _{\cdot }\left( t\pm \gamma \right)
\right\vert ^{p}dt\right] ^{1/p}\right\Vert _{S^{p}}\leq \omega f\left(
\gamma \right) _{S^{p}}
\end{equation*}%
and%
\begin{equation*}
\left\Vert w_{\cdot }f\left( \delta \right) _{p}\right\Vert _{S^{p}}\leq
\omega f\left( \delta \right) _{S^{p}},
\end{equation*}%
the analysis of the proof of Proposition 4 shows that, the estimate from
Theorem 5 implies the estimate from Theorem 2 with $p\geq q$ (without change 
$q$ instead of $q^{\prime }$\ in the estimate of the quantity $\left\{ \frac{%
1}{n+1}\sum_{k=n}^{2n}\left\vert I_{3}(k)\right\vert ^{q}\right\} ^{1/q}$).
Thus, taking $a_{n,k}=\frac{1}{n+1}$ when $k\leq n$ and $a_{n,k}=0$
otherwise, in the case $p=\infty $, we obtain the better estimate than this
one from Theorem 3 \cite[Theorem 5]{PT}.
\end{remark}

\section{Proofs of the results}

\subsection{Proof of Proposition 4}

In the proof we will use the following function $\Phi _{x}f\left( \delta
,\nu \right) =\frac{1}{\delta }\int_{\nu }^{\nu +\delta }\varphi _{x}\left(
u\right) du,$ with $\delta =\delta _{n}=\frac{\pi }{n+1}$ and its estimate
from \cite[Lemma 1, p.218]{WL}%
\begin{equation}
\left\vert \Phi _{x}f\left( \delta _{1},\delta _{2}\right) \right\vert \leq
w_{x}\left( \delta _{1}\right) +w_{x}\left( \delta _{2}\right)  \label{W}
\end{equation}%
for $f\in \Omega _{\alpha ,p}\left( w_{x}\right) $ and any $\delta
_{1},\delta _{2}>0$.

Since, for $n=0$ our estimate is evident we consider $n>0$, only.

Denote by $S_{k}^{\ast }f$ the sums of the form 
\begin{equation*}
S_{\frac{\alpha k}{2}}f\left( x\right) =\sum_{\left\vert \lambda _{\nu
}\right\vert \leq \frac{\alpha k}{2}}A_{\nu }\left( f\right) e^{i\lambda
_{\nu }x}
\end{equation*}%
such that the interval $\left( \frac{\alpha k}{2},\frac{\alpha \left(
k+1\right) }{2}\right) $ does not contain any $\lambda _{\nu }.$ Applying
Lemma 1.10.2 of \cite{BML} we easily verify that 
\begin{equation*}
S_{k}^{\ast }f\left( x\right) -f\left( x\right) =\int_{0}^{\infty }\varphi
_{x}\left( t\right) \Psi _{k}\left( t\right) dt,
\end{equation*}%
where\ $\varphi _{x}\left( t\right) :=f\left( x+t\right) +f\left( x-t\right)
-2f\left( x\right) $ and $\Psi _{k}\left( t\right) =\Psi _{\frac{\alpha k}{2}%
,\frac{\alpha \left( k+1\right) }{2}}\left( t\right) ,$ i.e. 
\begin{equation*}
\Psi _{k}\left( t\right) =\frac{4\sin \frac{\alpha t}{4}\sin \frac{\alpha
\left( 2k+1\right) t}{4}}{\alpha \pi t^{2}}
\end{equation*}%
(see also \cite{ASB}, p.41). Evidently, if the interval $\left( \frac{\alpha
k}{2},\frac{\alpha \left( k+1\right) }{2}\right) $ contains a Fourier
exponent $\lambda _{\nu },$ then 
\begin{equation*}
S_{\frac{\alpha k}{2}}f\left( x\right) =S_{k+1}^{\ast }f\left( x\right)
-\left( A_{\nu }\left( f\right) e^{i\lambda _{\nu }x}+A_{-\nu }\left(
f\right) e^{-i\lambda _{\nu }x}\right) .
\end{equation*}%
We can also note that if $f\in \Omega _{\alpha ,p}\left( w_{x}\right) ,$
with $p>1$ and $q>0$, then there exists $q^{\prime }\geq q$\ such that $%
q^{\prime }\geq 2$ and $p^{\prime }=\frac{q^{\prime }}{q^{\prime }-1}\leq p.$
Thus, 
\begin{equation*}
\left\{ \frac{1}{n+1}\sum_{k=n}^{2n}\left\vert S_{\frac{\alpha k}{2}}f\left(
x\right) -f\left( x\right) \right\vert ^{q}\right\} ^{1/q}\leq \left\{ \frac{%
1}{n+1}\sum_{k=n}^{2n}\left\vert S_{\frac{\alpha k}{2}}f\left( x\right)
-f\left( x\right) \right\vert ^{q^{\prime }}\right\} ^{1/q^{\prime }}.
\end{equation*}%
Since (see \cite[p.78]{ABI} and \cite[p. 7]{ADB}) 
\begin{equation*}
\left\{ \sum_{\nu =-\infty }^{\infty }\left\vert A_{\nu }\left( f\right)
\right\vert ^{q^{\prime }}\right\} ^{1/q^{\prime }}\leq \left\Vert
f\right\Vert _{B^{p^{\prime }}}\text{ \ \ and \ \ }\left\Vert f\right\Vert
_{B^{p^{\prime }}}\leq \left\Vert f\right\Vert _{S^{p^{\prime }}}\text{,}
\end{equation*}%
where $\left\Vert \cdot \right\Vert _{B^{p^{\prime },}}$\ with $(p^{\prime
}\geq 1),$ is the Besicovitch norm,\ so we have 
\begin{equation*}
\left\vert A_{\pm \nu }\left( f\right) \right\vert =\left\vert A_{\pm \nu
}\left( f-g_{\alpha \mu /2}\right) \right\vert \leq \left\Vert f-g_{\alpha
\mu /2}\right\Vert _{S^{p^{\prime }}}\leq \left\Vert f-g_{\alpha \mu
/2}\right\Vert _{S^{p}}=E_{\alpha \mu /2}\left( f\right) _{S^{p}},
\end{equation*}%
for some\ $g_{\alpha \mu /2}\in B_{\alpha \mu /2},$ with $\alpha k/2<\alpha
\mu /2<\lambda _{\nu }.$ Therefore, the deviation 
\begin{equation*}
\left\{ \frac{1}{n+1}\sum_{k=n}^{2n}\left\vert S_{\frac{\alpha k}{2}}f\left(
x\right) -f\left( x\right) \right\vert ^{q^{\prime }}\right\} ^{\frac{1}{%
q^{\prime }}}
\end{equation*}%
can be estimated from above by 
\begin{eqnarray*}
&&\left\{ \frac{1}{n+1}\sum_{k=n}^{2n}\left\vert \int_{0}^{\infty }\varphi
_{x}\left( t\right) \Psi _{k+\kappa }\left( t\right) dt\right\vert
^{q^{\prime }}\right\} ^{1/q^{\prime }}+\left\{ \frac{1}{n+1}%
\sum_{k=n}^{2n}\left( E_{\alpha k/2}\left( f\right) _{S^{p}}\right)
^{q^{\prime }}\right\} ^{1/q^{\prime }} \\
&\leq &\left\{ \frac{1}{n+1}\sum_{k=n}^{2n}\left\vert \int_{0}^{\infty
}\varphi _{x}\left( t\right) \Psi _{k+\kappa }\left( t\right) dt\right\vert
^{q^{\prime }}\right\} ^{1/q^{\prime }}+E_{\alpha n/2}\left( f\right)
_{S^{p}},
\end{eqnarray*}%
where $\kappa $ equals $0$ or $1.$ Putting $h_{n}=2\pi /\alpha (n+1)$ we
obtain 
\begin{eqnarray*}
&&\left\{ \frac{1}{n+1}\sum_{k=n}^{2n}\left\vert \int_{0}^{\infty }\varphi
_{x}\left( t\right) \Psi _{k+\kappa }\left( t\right) dt\right\vert
^{q^{\prime }}\right\} ^{1/q^{\prime }} \\
&=&\left\{ \frac{1}{n+1}\sum_{k=n}^{2n}\left\vert \left(
\int_{0}^{h_{n}}+\int_{h_{n}}^{(n+1)h_{n}}+\int_{(n+1)h_{n}}^{\infty
}\right) \varphi _{x}\left( t\right) \Psi _{k+\kappa }\left( t\right)
dt\right\vert ^{q^{\prime }}\right\} ^{1/q^{\prime }} \\
&\leq &\left\{ \frac{1}{n+1}\sum_{k=n}^{2n}\left\vert I_{1}(k)\right\vert
^{q^{\prime }}\right\} ^{1/q^{\prime }}+\left\{ \frac{1}{n+1}%
\sum_{k=n}^{2n}\left\vert I_{2}(k)\right\vert ^{q^{\prime }}\right\}
^{1/q^{\prime }}+\left\{ \frac{1}{n+1}\sum_{k=n}^{2n}\left\vert
I_{3}(k)\right\vert ^{q^{\prime }}\right\} ^{1/q^{\prime }}.
\end{eqnarray*}

So, for the first term we have 
\begin{equation*}
\left\{ \frac{1}{n+1}\sum_{k=n}^{2n}\left\vert I_{1}(k)\right\vert
^{q^{\prime }}\right\} ^{1/q^{\prime }}
\end{equation*}%
\begin{eqnarray*}
&\leq &\left\{ \frac{1}{n+1}\sum_{k=n}^{2n}\left\vert \frac{4}{\alpha \pi }%
\int_{0}^{h_{n}}\varphi _{x}\left( t\right) \frac{\sin \frac{\alpha t}{4}}{%
t^{2}}\sin \frac{\alpha t}{4}(2k+2\kappa +1)dt\right\vert ^{q^{\prime
}}\right\} ^{1/q^{\prime }} \\
&\leq &\frac{\alpha \left( 4n+2\kappa +1\right) }{4\pi }\int_{0}^{h_{n}}%
\left\vert \varphi _{x}\left( t\right) \right\vert dt\ll \left( 2+\frac{%
\kappa }{n+1}\right) w_{x}f\left( \frac{2\pi }{\alpha (n+1)}\right) _{p} \\
&\ll &w_{x}\left( \frac{\pi }{n+1}\right) .
\end{eqnarray*}%
Next, we estimate the second term. We have 
\begin{equation*}
\left\{ \frac{1}{n+1}\sum_{k=n}^{2n}\left\vert I_{2}(k)\right\vert
^{q^{\prime }}\right\} ^{1/q^{\prime }}
\end{equation*}%
\begin{eqnarray*}
&\leq &\left\{ \frac{1}{n+1}\sum_{k=n}^{2n}\left\vert \frac{4}{\alpha \pi }%
\int_{h_{n}}^{(n+1)h_{n}}\frac{\left( \Phi _{x}f\left( \delta _{n},t\right)
-\varphi _{x}\left( t\right) \right) \sin \frac{\alpha t}{4}}{t^{2}}\sin 
\frac{\alpha t}{4}(2k+2\kappa +1)dt\right\vert ^{q^{\prime }}\right\}
^{1/q^{\prime }} \\
&&+\left\{ \frac{1}{n+1}\sum_{k=n}^{2n}\left\vert \frac{4}{\alpha \pi }%
\int_{h_{n}}^{(n+1)h_{n}}\frac{\Phi _{x}f\left( \delta _{n},t\right) \sin 
\frac{\alpha t}{4}}{t^{2}}\sin \frac{\alpha t}{4}(2k+2\kappa
+1)dt\right\vert ^{q^{\prime }}\right\} ^{1/q^{\prime }} \\
&=&\left\{ \frac{1}{n+1}\sum_{k=n}^{2n}\left\vert I_{21}(k)\right\vert
^{q^{\prime }}\right\} ^{1/q^{\prime }}+\left\{ \frac{1}{n+1}%
\sum_{k=n}^{2n}\left\vert I_{22}(k)\right\vert ^{q^{\prime }}\right\}
^{1/q^{\prime }}.
\end{eqnarray*}%
From the Hausdorff -Young inequality \cite[Chap. XII, Th. 3.3 II]{Z} we
obtain 
\begin{equation*}
\left\{ \frac{1}{n+1}\sum_{k=n}^{2n}\left\vert I_{21}(k)\right\vert
^{q^{\prime }}\right\} ^{1/q^{\prime }}
\end{equation*}%
\begin{eqnarray*}
&\leq &\frac{8}{\alpha ^{2}}\left\{ \frac{1}{n+1}\sum_{k=n}^{2n}\left\vert 
\frac{\alpha }{2\pi }\int_{h_{n}}^{(n+1)h_{n}}\frac{\left( \Phi _{x}f\left(
\delta _{n},t\right) -\varphi _{x}\left( t\right) \right) \sin \frac{\alpha t%
}{4}}{t^{2}}\sin \frac{\alpha t}{4}(2\kappa +1)\cos \frac{\alpha kt}{2}%
dt\right\vert ^{q^{\prime }}\right\} ^{1/q^{\prime }} \\
&&+\frac{8}{\alpha ^{2}}\left\{ \frac{1}{n+1}\sum_{k=n}^{2n}\left\vert \frac{%
\alpha }{2\pi }\int_{h_{n}}^{(n+1)h_{n}}\frac{\left( \Phi _{x}f\left( \delta
_{n},t\right) -\varphi _{x}\left( t\right) \right) \sin \frac{\alpha t}{4}}{%
t^{2}}\cos \frac{\alpha t}{4}(2\kappa +1)\sin \frac{\alpha kt}{2}%
dt\right\vert ^{q^{\prime }}\right\} ^{1/q^{\prime }}
\end{eqnarray*}%
\begin{eqnarray*}
&\ll &\frac{1}{\left( n+1\right) ^{1/q^{\prime }}}\left( \left\{
\int_{h_{n}}^{(n+1)h_{n}}\left\vert \frac{\left( \Phi _{x}f\left( \delta
_{n},t\right) -\varphi _{x}\left( t\right) \right) \sin \frac{\alpha t}{4}}{%
t^{2}}\sin \frac{\alpha t}{4}(2\kappa +1)\right\vert ^{p^{\prime
}}dt\right\} ^{1/p^{\prime }}\right.  \\
&&+\left. \left\{ \int_{h_{n}}^{(n+1)h_{n}}\left\vert \frac{\left( \Phi
_{x}f\left( \delta _{n},t\right) -\varphi _{x}\left( t\right) \right) \sin 
\frac{\alpha t}{4}}{t^{2}}\cos \frac{\alpha t}{4}(2\kappa +1)\right\vert
^{p^{\prime }}dt\right\} ^{1/p^{\prime }}\right) 
\end{eqnarray*}

\begin{eqnarray*}
&\ll &\frac{1}{\left( n+1\right) ^{^{1/q^{\prime }}}}\left\{
\int_{h_{n}}^{(n+1)h_{n}}\left\vert \frac{\Phi _{x}f\left( \delta
_{n},t\right) -\varphi _{x}\left( t\right) }{t}\right\vert ^{p^{\prime
}}dt\right\} ^{1/p^{\prime }} \\
&=&\frac{1}{\left( n+1\right) ^{^{1/q^{\prime }}}}\left\{
\int_{h_{n}}^{(n+1)h_{n}}\left\vert \frac{1}{t\delta _{n}}\int_{t}^{t+\delta
n}\left( \varphi _{x}\left( u\right) -\varphi _{x}\left( t\right) \right)
du\right\vert ^{p^{\prime }}dt\right\} ^{1/p^{\prime }} \\
&\leq &\frac{1}{\left( n+1\right) ^{^{1/q^{\prime }}}}\left\{
\int_{h_{n}}^{(n+1)h_{n}}\left[ \frac{1}{t\delta _{n}}\int_{0}^{\delta
_{n}}\left\vert \varphi _{x}\left( u+t\right) -\varphi _{x}\left( t\right)
\right\vert du\right] ^{p^{\prime }}dt\right\} ^{1/p^{\prime }}
\end{eqnarray*}%
and by the Minkowski inequality, for $p^{\prime }>1$, we have 
\begin{equation*}
\left\{ \frac{1}{n+1}\sum_{k=n}^{2n}\left\vert I_{21}(k)\right\vert
^{q^{\prime }}\right\} ^{1/q^{\prime }}
\end{equation*}%
\begin{eqnarray*}
&\leq &\frac{1}{\left( n+1\right) ^{^{1/q^{\prime }}}}\frac{1}{\delta _{n}}%
\int_{0}^{\delta _{n}}\left[ \int_{h_{n}}^{(n+1)h_{n}}\frac{\left\vert
\varphi _{x}\left( u+t\right) -\varphi _{x}\left( t\right) \right\vert
^{p^{\prime }}}{t^{p^{\prime }}}dt\right] ^{1/p^{\prime }}du \\
&=&\frac{1}{\left( n+1\right) ^{^{1/q^{\prime }}}}\frac{1}{\delta _{n}}%
\int_{0}^{\delta _{n}}\left[ \int_{h_{n}}^{(n+1)h_{n}}\frac{1}{t^{p^{\prime
}}}\frac{d}{dt}\left( \int_{0}^{t}\left\vert \varphi _{x}\left( u+s\right)
-\varphi _{x}\left( s\right) \right\vert ^{p^{\prime }}ds\right) dt\right]
^{1/p^{\prime }}du
\end{eqnarray*}%
\begin{eqnarray*}
&=&\frac{1}{\left( n+1\right) ^{^{1/q^{\prime }}}}\frac{1}{\delta _{n}}%
\int_{0}^{\delta _{n}}\left\{ \left[ \frac{1}{t^{p^{\prime }}}%
\int_{0}^{t}\left\vert \varphi _{x}\left( u+s\right) -\varphi _{x}\left(
s\right) \right\vert ^{p^{\prime }}ds\right] _{h_{n}}^{(n+1)h_{n}}\right.  \\
&&+\left. p^{\prime }\int_{h_{n}}^{(n+1)h_{n}}\frac{1}{t^{p^{\prime }}}%
\left( \frac{1}{t}\int_{0}^{t}\left\vert \varphi _{x}\left( u+s\right)
-\varphi _{x}\left( s\right) \right\vert ^{p^{\prime }}ds\right) dt\right\}
^{1/p^{\prime }}du
\end{eqnarray*}%
\begin{eqnarray*}
&\ll &\frac{1}{\left( n+1\right) ^{^{1/q^{\prime }}}}\frac{1}{\delta _{n}}%
\int_{0}^{\delta _{n}}\left\{ \left( \frac{\alpha }{2\pi }\right)
^{p^{\prime }}\int_{0}^{(n+1)h_{n}}\left\vert \varphi _{x}\left( u+s\right)
-\varphi _{x}\left( s\right) \right\vert ^{p^{\prime }}ds\right. + \\
&&+\left( \frac{\alpha \left( n+1\right) }{2\pi }\right) ^{p^{\prime
}}\int_{0}^{h_{n}}\left\vert \varphi _{x}\left( u+s\right) -\varphi
_{x}\left( s\right) \right\vert ^{p^{\prime }}ds+ \\
&&+\left. p^{\prime }\int_{h_{n}}^{(n+1)h_{n}}\frac{1}{t^{p^{\prime }}}%
\left( w_{x}\left( u\right) \right) ^{p^{\prime }}dt\right\} ^{1/p^{\prime
}}du
\end{eqnarray*}%
\begin{eqnarray*}
&\ll &\frac{1}{\left( n+1\right) ^{^{1/q^{\prime }}}}\frac{1}{\delta _{n}}%
\int_{0}^{\delta _{n}}\left\{ w_{x}^{p^{\prime }}\left( u\right) \left(
\left( \frac{\alpha }{2\pi }\right) ^{p^{\prime }-1}+\left( \frac{\alpha
\left( n+1\right) }{2\pi }\right) ^{p^{\prime }-1}\right. \right.  \\
&&+\left. \left. p^{\prime }\int_{h_{n}}^{(n+1)h_{n}}\frac{1}{t^{p^{\prime }}%
}dt\right) \right\} ^{1/p^{\prime }}du \\
&\ll &\frac{1}{\left( n+1\right) ^{^{1/q^{\prime }}}}w_{x}\left( \delta
_{n}\right) \left\{ 1+\left( n+1\right) ^{p^{\prime }-1}+p^{\prime }\left[ 
\frac{t^{1-p^{\prime }}}{1-p^{\prime }}\right] _{h_{n}}^{(n+1)h_{n}}\right\}
^{1/p^{\prime }} \\
&\ll &w_{x}\left( \delta _{n}\right) \left( n+1\right) ^{-1/q^{\prime
}}\left( n+1\right) ^{1-1/p^{\prime }}=w_{x}\left( \frac{\pi }{n+1}\right) .
\end{eqnarray*}

And 
\begin{equation*}
\left\{ \frac{1}{n+1}\sum_{k=n}^{2n}\left\vert I_{22}(k)\right\vert
^{q^{\prime }}\right\} ^{1/q^{\prime }}
\end{equation*}%
\begin{eqnarray*}
&=&\left\{ \frac{1}{n+1}\sum_{k=n}^{2n}\left\vert \frac{4}{\alpha \pi }%
\int_{h_{n}}^{(n+1)h_{n}}\frac{\Phi _{x}f\left( \delta _{n},t\right) \sin 
\frac{\alpha t}{4}}{t^{2}}\sin \frac{\alpha t}{4}(2k+2\kappa
+1)dt\right\vert ^{q^{\prime }}\right\} ^{1/q^{\prime }} \\
&=&\left\{ \frac{1}{n+1}\sum_{k=n}^{2n}\left\vert \frac{4}{\alpha \pi }%
\int_{h_{n}}^{(n+1)h_{n}}\Phi _{x}f\left( \delta _{n},t\right) \frac{\sin 
\frac{\alpha t}{4}}{t^{2}}\frac{d}{dt}\left( \frac{-\cos \frac{\alpha t}{4}%
(2k+2\kappa +1)}{\frac{\alpha }{4}(2k+2\kappa +1)}\right) dt\right\vert
^{q^{\prime }}\right\} ^{1/q^{\prime }}
\end{eqnarray*}%
\begin{eqnarray*}
&=&\left\{ \frac{1}{n+1}\sum_{k=n}^{2n}\left\vert \frac{4}{\alpha \pi }\left[
\Phi _{x}f\left( \delta _{n},t\right) \frac{\sin \frac{\alpha t}{4}}{t^{2}}%
\left( \frac{-\cos \frac{\alpha t}{4}(2k+2\kappa +1)}{\frac{\alpha }{4}%
(2k+2\kappa +1)}\right) \right] _{h_{n}}^{(n+1)h_{n}}\right. \right.  \\
&&\left. \left. +\frac{4}{\alpha \pi }\int_{h_{n}}^{(n+1)h_{n}}\frac{d}{dt}%
\left( \Phi _{x}f\left( \delta _{n},t\right) \frac{\sin \frac{\alpha t}{4}}{%
t^{2}}\right) \frac{\cos \frac{\alpha t}{4}(2k+2\kappa +1)}{\frac{\alpha }{4}%
(2k+2\kappa +1)}dt\right\vert ^{q^{\prime }}\right\} ^{1/q^{\prime }}
\end{eqnarray*}%
\begin{eqnarray*}
&\leq &\left\{ \frac{1}{n+1}\sum_{k=n}^{2n}\left\vert \frac{-4}{\alpha \pi }%
\Phi _{x}f\left( \delta _{n},\frac{2\pi }{\alpha }\right) \left( \frac{%
\alpha }{2\pi }\right) ^{2}\frac{\cos \frac{\pi }{2}(2k+2\kappa +1)}{\frac{%
\alpha }{4}(2k+2\kappa +1)}\right. \right.  \\
&&+\frac{4}{\alpha \pi }\left. \left. \Phi _{x}f\left( \delta _{n},\frac{%
2\pi }{\alpha (n+1)}\right) \frac{\sin \frac{\pi }{2(n+1)}}{\left( \frac{%
2\pi }{\alpha (n+1)}\right) ^{2}}\frac{\cos \frac{\pi }{2(n+1)}(2k+2\kappa
+1)}{\frac{\alpha }{4}(2k+2\kappa +1)}\right\vert ^{q^{\prime }}\right\}
^{1/q^{\prime }} \\
&&+\left\{ \frac{1}{n+1}\sum_{k=n}^{2n}\left\vert \frac{4}{\alpha \pi }%
\int_{h_{n}}^{(n+1)h_{n}}\frac{d}{dt}\left( \Phi _{x}f\left( \delta
_{n},t\right) \frac{\sin \frac{\alpha t}{4}}{t^{2}}\right) \frac{\cos \frac{%
\alpha t}{4}(2k+2\kappa +1)}{\frac{\alpha }{4}(2k+2\kappa +1)}dt\right\vert
^{q^{\prime }}\right\} ^{1/q^{\prime }}
\end{eqnarray*}%
\begin{equation*}
=\left\{ \frac{1}{n+1}\sum_{k=n}^{2n}\left\vert I_{221}(k)\right\vert
^{q^{\prime }}\right\} ^{1/q^{\prime }}+\left\{ \frac{1}{n+1}%
\sum_{k=n}^{2n}\left\vert I_{222}(k)\right\vert ^{q^{\prime }}\right\}
^{1/q^{\prime }}.
\end{equation*}%
For the first term, using (\ref{W}), we obtain 
\begin{equation*}
\left\{ \frac{1}{n+1}\sum_{k=n}^{2n}\left\vert I_{221}(k)\right\vert
^{q^{\prime }}\right\} ^{1/q^{\prime }}
\end{equation*}%
\begin{eqnarray*}
&\leq &\left\{ \frac{1}{n+1}\sum_{k=n}^{2n}\left( \frac{4}{\pi ^{3}}%
\left\vert \Phi _{x}f\left( \delta _{n},\frac{2\pi }{\alpha }\right)
\right\vert \frac{1}{(2k+2\kappa +1)}\right) ^{q^{\prime }}\right\}
^{1/q^{\prime }} \\
&&+\left\{ \frac{1}{n+1}\sum_{k=n}^{2n}\left( \frac{2(n+1)}{\pi ^{2}}%
\left\vert \Phi _{x}f\left( \delta _{n},\frac{2\pi }{\alpha (n+1)}\right)
\right\vert \frac{1}{(2k+2\kappa +1)}\right) ^{q^{\prime }}\right\}
^{1/q^{\prime }}
\end{eqnarray*}%
\begin{eqnarray*}
&\leq &\left\{ \frac{1}{n+1}\sum_{k=n}^{2n}\left( \frac{4}{\left( k+1\right)
\pi ^{3}}\left( w_{x}\left( \delta _{n}\right) +w_{x}\left( \frac{2\pi }{%
\alpha }\right) \right) \right) ^{q^{\prime }}\right\} ^{1/q^{\prime }} \\
&&\left\{ \frac{1}{n+1}\sum_{k=n}^{2n}\left( \frac{2(n+1)}{\left( k+1\right)
\pi ^{2}}\left( w_{x}\left( \delta _{n}\right) +w_{x}\left( \frac{2\pi }{%
\alpha (n+1)}\right) \right) \right) ^{q^{\prime }}\right\} ^{1/q^{\prime }}
\end{eqnarray*}%
\begin{eqnarray*}
&\ll &w_{x}\left( \delta _{n}\right) +w_{x}\left( \frac{2\pi }{\alpha (n+1)}%
\right) +\left\{ \frac{1}{n+1}\sum_{k=n}^{2n}\left( \frac{1}{k+1}w_{x}\left( 
\frac{2\pi }{\alpha }\right) \right) ^{q^{\prime }}\right\} ^{1/q^{\prime }}
\\
&\ll &w_{x}\left( \frac{\pi }{n+1}\right) +\left\{ \frac{1}{n+1}%
\sum_{k=n}^{2n}w_{x}^{q^{\prime }}\left( \frac{\pi }{k+1}\right) \right\}
^{1/q^{\prime }}\leq 2w_{x}\left( \frac{\pi }{n+1}\right) .
\end{eqnarray*}

For the second term, using (\ref{W}) and the Hausdorff-Young inequality \cite%
[Chap. XII, Th. 3.3 II]{Z} we have 
\begin{equation*}
\left\{ \frac{1}{n+1}\sum_{k=n}^{2n}\left\vert I_{222}(k)\right\vert
^{q^{\prime }}\right\} ^{1/q^{\prime }}
\end{equation*}%
\begin{eqnarray*}
&\leq &\frac{1}{\left( n+1\right) ^{1/q^{\prime }+1}}\left\{
\sum_{k=n}^{2n}\left\vert \frac{16}{\alpha ^{2}\pi }\int_{h_{n}}^{(n+1)h_{n}}%
\frac{d}{dt}\left( \Phi _{x}f\left( \delta _{n},t\right) \frac{\sin \frac{%
\alpha t}{4}}{t^{2}}\right) \cos \frac{\alpha t}{4}(2k+2\kappa
+1)dt\right\vert ^{q^{\prime }}\right\} ^{1/q^{\prime }} \\
&\ll &\frac{1}{\left( n+1\right) ^{1/q^{\prime }+1}}\left\{
\int_{h_{n}}^{(n+1)h_{n}}\left\vert \frac{d}{dt}\left( \Phi _{x}f\left(
\delta _{n},t\right) \frac{\sin \frac{\alpha t}{4}}{t^{2}}\right)
\right\vert ^{p^{\prime }}dt\right\} ^{^{1/p^{\prime }}}
\end{eqnarray*}%
\begin{eqnarray*}
&\leq &\frac{1}{\left( n+1\right) ^{1/q^{\prime }+1}}\left\{
\int_{h_{n}}^{(n+1)h_{n}}\left( \frac{1}{\delta _{n}}\left\vert \varphi
_{x}\left( \delta _{n}+t\right) -\varphi _{x}\left( t\right) \right\vert 
\frac{\sin \frac{\alpha t}{4}}{t^{2}}\right. \right.  \\
&&\left. \left. +\left\vert \Phi _{x}f\left( \delta _{n},t\right)
\right\vert \frac{\left\vert \frac{\alpha }{4}t^{2}\cos \frac{\alpha t}{4}%
-2t\sin \frac{\alpha t}{4}\right\vert }{t^{4}}\right) ^{p^{\prime
}}dt\right\} ^{^{1/p^{\prime }}}
\end{eqnarray*}%
\begin{eqnarray*}
&\ll &\frac{1}{\left( n+1\right) ^{1/q^{\prime }+1}}\left[ \left\{
\int_{h_{n}}^{(n+1)h_{n}}\left( \frac{1}{\delta _{n}t}\left\vert \varphi
_{x}\left( \delta _{n}+t\right) -\varphi _{x}\left( t\right) \right\vert
\right) ^{p^{\prime }}dt\right\} ^{^{1/p^{\prime }}}\right.  \\
&&\left. +\left\{ \int_{h_{n}}^{(n+1)h_{n}}\left( \frac{1}{t^{2}}\left(
w_{x}\left( \delta _{n}\right) +w_{x}\left( t\right) \right) \right)
^{p^{\prime }}dt\right\} ^{^{1/p^{\prime }}}\right] 
\end{eqnarray*}%
\begin{eqnarray*}
&\leq &\frac{1}{\left( n+1\right) ^{1/q^{\prime }+1}}\left[ \left\{ \frac{1}{%
\delta _{n}^{p^{\prime }}}\left[ t^{-p^{\prime }}\int_{0}^{t}\left\vert
\varphi _{x}\left( \delta _{n}+u\right) -\varphi _{x}\left( u\right)
\right\vert ^{p^{\prime }}du\right] _{h_{n}}^{(n+1)h_{n}}\right. \right.  \\
&&\left. +\frac{p^{\prime }}{\delta _{n}^{p^{\prime }}}%
\int_{h_{n}}^{(n+1)h_{n}}t^{-p^{\prime }}\left( \frac{1}{t}%
\int_{0}^{t}\left\vert \varphi _{x}\left( \delta _{n}+u\right) -\varphi
_{x}\left( u\right) \right\vert ^{p^{\prime }}du\right) dt\right\}
^{^{1/p^{\prime }}} \\
&&+\left. \left\{ \int_{h_{n}}^{(n+1)h_{n}}\left( \frac{w_{x}\left( \delta
_{n}\right) }{t^{2}}\right) ^{p^{\prime }}dt\right\} ^{^{^{1/p^{\prime
}}}}+\left\{ \int_{h_{n}}^{(n+1)h_{n}}\left( \frac{w_{x}\left( t\right) /t}{t%
}\right) ^{p^{\prime }}dt\right\} ^{^{1/p^{\prime }}}\right]  \\
&&+\left. w_{x}\left( \delta _{n}\right) \left\{ \left[ \frac{t^{-2p^{\prime
}+1}}{-2p^{\prime }+1}\right] _{h_{n}}^{(n+1)h_{n}}\right\} ^{^{1/p^{\prime
}}}+(n+1)w_{x}\left( \frac{\pi }{n+1}\right) \left\{ \left[ \frac{%
t^{-p^{\prime }+1}}{-p^{\prime }+1}\right] _{h_{n}}^{(n+1)h_{n}}\right\}
^{^{1/p^{\prime }}}\right] 
\end{eqnarray*}%
\begin{eqnarray*}
&\ll &\frac{1}{\left( n+1\right) ^{1/q^{\prime }+1}}\left[ \frac{1}{\delta
_{n}}\left\{ \frac{\alpha }{2\pi }\int_{0}^{(n+1)h_{n}}\left\vert \varphi
_{x}\left( \delta _{n}+u\right) -\varphi _{x}\left( u\right) \right\vert
^{p^{\prime }}du\right\} ^{^{1/p^{\prime }}}\right.  \\
&&+\frac{1}{\delta _{n}}(n+1)^{1-1/p^{\prime }}\left\{ \frac{\alpha (n+1)}{%
2\pi }\int_{0}^{h_{n}}\left\vert \varphi _{x}\left( \delta _{n}+u\right)
-\varphi _{x}\left( u\right) \right\vert ^{p^{\prime }}du\right\}
^{^{1/p^{\prime }}} \\
&&+\frac{1}{\delta _{n}}\left\{ \int_{h_{n}}^{(n+1)h_{n}}t^{-p^{\prime
}}\left( w_{x}\left( \delta _{n}\right) \right) ^{p^{\prime }}dt\right\}
^{^{1/p^{\prime }}}
\end{eqnarray*}%
\begin{eqnarray*}
&\ll &\frac{1}{\left( n+1\right) ^{1/q^{\prime }+1}}\left[ \frac{1}{\delta
_{n}}w_{x}\left( \delta _{n}\right) +\frac{1}{\delta _{n}}w_{x}\left( \delta
_{n}\right) \left( n+1\right) ^{1-1/p^{\prime }}\right.  \\
&&\left. +w_{x}\left( \delta _{n}\right) \left( n+1\right) ^{2-1/p^{\prime
}}+(n+1)w_{x}\left( \frac{\pi }{n+1}\right) \left( n+1\right)
^{1-1/p^{\prime }}\right]  \\
&\ll &w_{x}\left( \frac{\pi }{n+1}\right) .
\end{eqnarray*}

For the third term we obtain 
\begin{equation*}
\left\{ \frac{1}{n+1}\sum_{k=n}^{2n}\left\vert I_{3}(k)\right\vert
^{q^{\prime }}\right\} ^{1/q^{\prime }}\leq 
\end{equation*}%
\begin{eqnarray*}
&\leq &\left\{ \frac{1}{n+1}\sum_{k=n}^{2n}\left\vert \sum_{\mu =1}^{\infty
}\int_{(n+1)h_{n}\mu }^{(n+1)h_{n}\left( \mu +1\right) }\left[ \varphi
_{x}\left( t\right) -\Phi _{x}f\left( \delta _{k},t\right) \right] \Psi
_{k+\kappa }\left( t\right) dt\right\vert ^{q^{\prime }}\right\}
^{1/q^{\prime }} \\
&&+\left\{ \frac{1}{n+1}\sum_{k=n}^{2n}\left\vert \sum_{\mu =1}^{\infty
}\int_{(n+1)h_{n}\mu }^{(n+1)h_{n}\left( \mu +1\right) }\Phi _{x}f\left(
\delta _{k},t\right) \Psi _{k+\kappa }\left( t\right) dt\right\vert
^{q^{\prime }}\right\} ^{1/q^{\prime }} \\
&=&\left\{ \frac{1}{n+1}\sum_{k=n}^{2n}\left\vert I_{31}(k)\right\vert
^{q^{\prime }}\right\} ^{1/q^{\prime }}+\left\{ \frac{1}{n+1}%
\sum_{k=n}^{2n}\left\vert I_{32}(k)\right\vert ^{q^{\prime }}\right\}
^{1/q^{\prime }}
\end{eqnarray*}%
and 
\begin{eqnarray*}
\left\vert I_{31}(k)\right\vert  &\leq &\frac{4}{\alpha \pi }\sum_{\mu
=1}^{\infty }\int_{(n+1)h_{n}\mu }^{(n+1)h_{n}\left( \mu +1\right)
}\left\vert \varphi _{x}\left( t\right) -\Phi _{x}f\left( \delta
_{k},t\right) \right\vert t^{-2}dt \\
&\leq &\frac{4}{\alpha \pi }\sum_{\mu =1}^{\infty }\int_{(n+1)h_{n}\mu
}^{(n+1)h_{n}\left( \mu +1\right) }\left[ \frac{1}{\delta _{k}t^{2}}%
\int_{0}^{\delta _{k}}\left\vert \varphi _{x}\left( t\right) -\varphi
_{x}\left( t+u\right) \right\vert du\right] dt
\end{eqnarray*}%
\begin{eqnarray*}
&=&\frac{4}{\alpha \pi }\frac{1}{\delta _{k}}\int_{0}^{\delta _{k}}\sum_{\mu
=1}^{\infty }\left\{ \int_{(n+1)h_{n}\mu }^{(n+1)h_{n}\left( \mu +1\right) }%
\frac{1}{t^{2}}\left\vert \varphi _{x}\left( t\right) -\varphi _{x}\left(
t+u\right) \right\vert dt\right\} du \\
&=&\frac{4}{\alpha \pi }\frac{1}{\delta _{k}}\int_{0}^{\delta _{k}}\sum_{\mu
=1}^{\infty }\left\{ \left[ \frac{1}{t^{2}}\int_{0}^{t}\left\vert \varphi
_{x}\left( s\right) -\varphi _{x}\left( s+u\right) \right\vert ds\right]
_{(n+1)h_{n}\mu }^{(n+1)h_{n}\left( \mu +1\right) }\right.  \\
&&+\left. 2\int_{(n+1)h_{n}\mu }^{(n+1)h_{n}\left( \mu +1\right) }\left[ 
\frac{1}{t^{3}}\int_{0}^{t}\left\vert \varphi _{x}\left( s\right) -\varphi
_{x}\left( s+u\right) \right\vert ds\right] dt\right\} du
\end{eqnarray*}%
\begin{eqnarray*}
&\ll &\left\vert \frac{1}{\delta _{k}}\int_{0}^{\delta _{k}}\sum_{\mu
=1}^{\infty }\left\{ \frac{1}{[(n+1)h_{n}\left( \mu +1\right) ]^{2}}%
\int_{0}^{(n+1)h_{n}\left( \mu +1\right) }\left\vert \varphi _{x}\left(
s\right) -\varphi _{x}\left( s+u\right) \right\vert ds\right. \right.  \\
&&\left. -\left. \frac{1}{[(n+1)h_{n}\mu ]^{2}}\int_{0}^{(n+1)h_{n}\mu
}\left\vert \varphi _{x}\left( s\right) -\varphi _{x}\left( s+u\right)
\right\vert ds\right\} du\right\vert  \\
&&+\frac{1}{\delta _{k}}\int_{0}^{\delta _{k}}\sum_{\mu =1}^{\infty }\left\{
\int_{(n+1)h_{n}\mu }^{(n+1)h_{n}\left( \mu +1\right) }\left[ \frac{1}{t^{3}}%
\int_{0}^{t}\left\vert \varphi _{x}\left( s\right) -\varphi _{x}\left(
s+u\right) \right\vert ds\right] dt\right\} du.
\end{eqnarray*}%
Since $f\in \Omega _{\alpha ,p}\left( w_{x}\right) $, thus for any $x$ 
\begin{equation*}
\lim_{\zeta \rightarrow \infty }\frac{1}{\zeta ^{2}}\int_{0}^{\zeta
}\left\vert \varphi _{x}\left( s\right) -\varphi _{x}\left( s+u\right)
\right\vert ds\leq \lim_{\zeta \rightarrow \infty }\frac{1}{\zeta }%
w_{x}\left( u\right) \leq \lim_{\zeta \rightarrow \infty }\frac{1}{\zeta }%
w_{x}\left( \delta _{k}\right) \leq \lim_{\zeta \rightarrow \infty }\frac{1}{%
\zeta }w_{x}\left( \pi \right) =0,
\end{equation*}%
and therefore 
\begin{eqnarray*}
\left\vert I_{31}(k)\right\vert  &\leq &\frac{1}{\delta _{k}}%
\int_{0}^{\delta _{k}}\frac{\alpha }{2\pi }\left[ \frac{\alpha }{2\pi }%
\int_{0}^{2\pi /\alpha }\left\vert \varphi _{x}\left( s\right) -\varphi
_{x}\left( s+u\right) \right\vert ds\right] du \\
&&+\frac{1}{\delta _{k}}\int_{0}^{\delta _{k}}w_{x}\left( u\right)
du\sum_{\mu =1}^{\infty }\left\{ \int_{(n+1)h_{n}\mu }^{(n+1)h_{n}\left( \mu
+1\right) }\frac{1}{t^{2}}dt\right\}  \\
&\ll &\frac{1}{\delta _{k}}\int_{0}^{\delta _{k}}w_{x}\left( u\right)
du+w_{x}\left( \delta _{k}\right) \sum_{\mu =1}^{\infty }\frac{1}{%
(n+1)h_{n}\mu ^{2}} \\
&\ll &w_{x}\left( \delta _{k}\right) .
\end{eqnarray*}%
Next, we will estimate the term $\left\vert I_{32}(k)\right\vert .$ So, 
\begin{eqnarray*}
I_{32}(k) &=&\frac{2}{\alpha \pi }\sum_{\mu =1}^{\infty }\int_{(n+1)h_{n}\mu
}^{(n+1)h_{n}\left( \mu +1\right) }\frac{\Phi _{x}f\left( \delta
_{k},t\right) }{t^{2}}\frac{d}{dt}\left( -\frac{\cos \frac{\alpha t\left(
k+\kappa \right) }{2}}{\frac{\alpha \left( k+\kappa \right) }{2}}+\frac{\cos 
\frac{\alpha t\left( k+\kappa +1\right) }{2}}{\frac{\alpha \left( k+\kappa
+1\right) }{2}}\right) dt \\
&=&\frac{2}{\alpha \pi }\sum_{\mu =1}^{\infty }\left[ \frac{\Phi _{x}f\left(
\delta _{k},t\right) }{t^{2}}\left( -\frac{\cos \frac{\alpha t\left(
k+\kappa \right) }{2}}{\frac{\alpha \left( k+\kappa \right) }{2}}+\frac{\cos 
\frac{\alpha t\left( k+\kappa +1\right) }{2}}{\frac{\alpha \left( k+\kappa
+1\right) }{2}}\right) \right] _{(n+1)h_{n}\mu }^{(n+1)h_{n}\left( \mu
+1\right) } \\
&&+\frac{2}{\alpha \pi }\sum_{\mu =1}^{\infty }\int_{(n+1)h_{n}\mu
}^{(n+1)h_{n}\left( \mu +1\right) }\frac{d}{dt}\left( \frac{\Phi _{x}f\left(
\delta _{k},t\right) }{t^{2}}\right) \left( \frac{\cos \frac{\alpha t\left(
k+\kappa \right) }{2}}{\frac{\alpha \left( k+\kappa \right) }{2}}-\frac{\cos 
\frac{\alpha t\left( k+\kappa +1\right) }{2}}{\frac{\alpha \left( k+\kappa
+1\right) }{2}}\right) dt \\
&=&I_{321}\left( k\right) +I_{322}\left( k\right) 
\end{eqnarray*}%
Since $f\in \Omega _{\alpha ,p}$, thus for $x$ (using (\ref{W}))%
\begin{eqnarray*}
&&\lim_{\zeta \rightarrow \infty }\left\vert \frac{\Phi _{x}f\left( \delta
_{k},\frac{2\pi }{\alpha }\zeta \right) }{\left[ \frac{2\pi }{\alpha }\zeta %
\right] ^{2}}\left( -\frac{\cos \left[ \pi \zeta (k+\kappa )\right] }{\frac{%
\alpha \left( k+\kappa \right) }{2}}+\frac{\cos \left[ \pi \zeta \left(
k+\kappa +1\right) \right] }{\frac{\alpha \left( k+\kappa +1\right) }{2}}%
\right) \right\vert  \\
&\leq &\lim_{\zeta \rightarrow \infty }\frac{w_{x}\left( \delta _{k}\right)
+w_{x}\left( \frac{2\pi }{\alpha }\zeta \right) }{2\pi ^{2}\zeta ^{2}k}\ll
\lim_{\zeta \rightarrow \infty }\frac{w_{x}\left( \delta _{k}\right) +\zeta
w_{x}\left( \frac{2\pi }{\alpha }\right) }{\zeta ^{2}k}\ll w_{x}\left( \pi
\right) \lim_{\zeta \rightarrow \infty }\frac{1+\zeta }{\zeta ^{2}}=0,
\end{eqnarray*}%
and therefore%
\begin{eqnarray*}
I_{321}\left( k\right)  &=&\frac{2}{\alpha \pi }\sum_{\mu =1}^{\infty }\left[
\frac{\Phi _{x}f\left( \delta _{k},\frac{2\pi }{\alpha }\left( \mu +1\right)
\right) }{\left[ \frac{2\pi }{\alpha }\left( \mu +1\right) \right] ^{2}}%
\left( -\frac{\cos \left[ \pi \left( \mu +1\right) \left( k+\kappa \right) %
\right] }{\frac{\alpha \left( k+\kappa \right) }{2}}\right. \right.  \\
&&+\left. \frac{\cos \left[ \pi \left( \mu +1\right) \left( k+\kappa
+1\right) \right] }{\frac{\alpha \left( k+\kappa +1\right) }{2}}\right)  \\
&&-\left. \frac{\Phi _{x}f\left( \delta _{k},\frac{2\pi }{\alpha }\mu
\right) }{\left[ \frac{2\pi }{\alpha }\mu \right] ^{2}}\left( -\frac{\cos %
\left[ \pi \mu (k+\kappa )\right] }{\frac{\alpha \left( k+\kappa \right) }{2}%
}+\frac{\cos \left[ \pi \mu \left( k+\kappa +1\right) \right] }{\frac{\alpha
\left( k+\kappa +1\right) }{2}}\right) \right] 
\end{eqnarray*}%
\begin{eqnarray*}
&=&-\frac{2}{\alpha \pi }\frac{\Phi _{x}f\left( \delta _{k},2\pi /\alpha
\right) }{\left[ 2\pi /\alpha \right] ^{2}}\left( -\frac{\left( -1\right)
^{\left( k+\kappa \right) }}{\frac{\alpha \left( k+\kappa \right) }{2}}+%
\frac{\left( -1\right) ^{\left( k+\kappa +1\right) }}{\frac{\alpha \left(
k+\kappa +1\right) }{2}}\right)  \\
&=&-\frac{1}{\pi ^{3}}\Phi _{x}f\left( \delta _{k},2\pi /\alpha \right)
\left( -1\right) ^{\left( k+\kappa +1\right) }\left( \frac{1}{k+\kappa +1}+%
\frac{1}{k+\kappa }\right) .
\end{eqnarray*}%
Using (\ref{W}), we get 
\begin{equation*}
\left\vert I_{321}\left( k\right) \right\vert \ll \frac{1}{\pi ^{3}}\frac{2}{%
k+1}\left\vert \Phi _{x}f\left( \delta _{k},2\pi /\alpha \right) \right\vert
\leq \frac{2}{\pi ^{3}\left( k+1\right) }\left( w_{x}\left( \delta
_{k}\right) +w_{x}\left( 2\pi /\alpha \right) \right) .
\end{equation*}%
Similarly 
\begin{eqnarray*}
I_{322}\left( k\right)  &=&\frac{2}{\alpha \pi }\sum_{\mu =1}^{\infty
}\int_{(n+1)h_{n}\mu }^{(n+1)h_{n}\left( \mu +1\right) }\left( \frac{\frac{d%
}{dt}\Phi _{x}f\left( \delta _{k},t\right) }{t^{2}}-\frac{2\Phi _{x}f\left(
\delta _{k},t\right) }{t^{3}}\right)  \\
&&\cdot \left( \frac{\cos \frac{\alpha t\left( k+\kappa \right) }{2}}{\frac{%
\alpha \left( k+\kappa \right) }{2}}-\frac{\cos \frac{\alpha t\left(
k+\kappa +1\right) }{2}}{\frac{\alpha \left( k+\kappa +1\right) }{2}}\right)
dt
\end{eqnarray*}%
and 
\begin{eqnarray*}
\left\vert I_{322}\left( k\right) \right\vert  &\ll &\frac{8}{\alpha
^{2}\left( k+1\right) \pi }\sum_{\mu =1}^{\infty }\left[ \int_{(n+1)h_{n}\mu
}^{(n+1)h_{n}\left( \mu +1\right) }\frac{\left\vert \varphi _{x}\left(
t+\delta _{k}\right) -\varphi _{x}\left( t\right) \right\vert }{\delta
_{k}t^{2}}dt\right.  \\
&&+\left. 2\int_{(n+1)h_{n}\mu }^{(n+1)h_{n}\left( \mu +1\right) }\frac{%
\left\vert \Phi _{x}f\left( \delta _{k},t\right) \right\vert }{t^{3}}dt%
\right]  \\
&\leq &\frac{8}{\alpha ^{2}\left( k+1\right) \pi \delta _{k}}\sum_{\mu
=1}^{\infty }\int_{(n+1)h_{n}\mu }^{(n+1)h_{n}\left( \mu +1\right) }\frac{%
\left\vert \varphi _{x}\left( t+\delta _{k}\right) -\varphi _{x}\left(
t\right) \right\vert }{t^{2}}dt \\
&&+\frac{16}{\alpha ^{2}\left( k+1\right) \pi }\sum_{\mu =1}^{\infty
}\int_{(n+1)h_{n}\mu }^{(n+1)h_{n}\left( \mu +1\right) }\frac{w_{x}\left(
\delta _{k}\right) +w_{x}\left( t\right) }{t^{3}}dt
\end{eqnarray*}%
\begin{eqnarray*}
&\ll &\frac{1}{\left( k+1\right) \delta _{k}}w_{x}\left( \delta _{k}\right) +%
\frac{1}{k+1}\sum_{\mu =1}^{\infty }\left[ \left( w_{x}\left( \delta
_{k}\right) +w_{x}\left( \frac{2\pi \left( \mu +1\right) }{\alpha }\right)
\right) \frac{\alpha ^{2}}{4\pi ^{2}\mu ^{3}}\right]  \\
&\ll &w_{x}\left( \delta _{k}\right) +\frac{1}{k+1}\left[ w_{x}\left( \delta
_{k}\right) \sum_{\mu =1}^{\infty }\frac{1}{\mu ^{3}}+\sum_{\mu =1}^{\infty }%
\frac{w_{x}\left( \frac{2\pi \left( \mu +1\right) }{\alpha }\right) }{\mu
^{3}}\right]  \\
&\ll &w_{x}\left( \delta _{k}\right) +\frac{1}{k+1}\left( w_{x}\left( \delta
_{k}\right) +w_{x}\left( \frac{4\pi }{\alpha }\right) \sum_{\mu =1}^{\infty }%
\frac{\mu +1}{\mu ^{3}}\right)  \\
&\ll &w_{x}\left( \delta _{k}\right) +\frac{1}{k+1}\left( w_{x}\left( \delta
_{k}\right) +w_{x}\left( \frac{4\pi }{\alpha }\right) \right) .
\end{eqnarray*}%
Therefore 
\begin{equation*}
\left\vert I_{3}\left( k\right) \right\vert \ll w_{x}\left( \delta
_{k}\right) +\frac{1}{k+1}\left( w_{x}\left( \delta _{k}\right) +w_{x}\left( 
\frac{2\pi }{\alpha }\right) +w_{x}\left( \frac{4\pi }{\alpha }\right)
\right) 
\end{equation*}%
and thus 
\begin{eqnarray*}
\left\{ \frac{1}{n+1}\sum_{k=n}^{2n}\left\vert I_{3}\left( k\right)
\right\vert ^{q^{\prime }}\right\} ^{1/q^{\prime }} &\ll &\left\{ \frac{1}{%
n+1}\sum_{k=n}^{2n}\left( w_{x}\left( \frac{\pi }{k+1}\right) +\frac{1}{k+1}%
w_{x}\left( \frac{\pi }{\alpha }\right) \right) ^{q^{\prime }}\right\}
^{1/q^{\prime }} \\
&\ll &\left\{ \frac{1}{n+1}\sum_{k=n}^{2n}\left( w_{x}\left( \frac{\pi }{k+1}%
\right) \right) ^{q^{\prime }}\right\} ^{1/q^{\prime }}\leq w_{x}\left( 
\frac{\pi }{n+1}\right) .
\end{eqnarray*}

and the desired result follows. $\square $

\subsection{Proof of Theorem 5}

For some $c>1$%
\begin{equation*}
H_{n,A,\gamma }^{q}f\left( x\right) =\left\{ \sum_{k=0}^{2^{\left[ c\right]
}-1}a_{n,k}\left\vert S_{\frac{\alpha k}{2}}f\left( x\right) -f\left(
x\right) \right\vert ^{q}+\sum_{k=2^{[c]}}^{\infty }a_{n,k}\left\vert S_{%
\frac{\alpha k}{2}}f\left( x\right) -f\left( x\right) \right\vert
^{q}\right\} ^{1/q}
\end{equation*}%
\begin{eqnarray*}
&\ll &\left\{ \sum_{k=0}^{2^{\left[ c\right] }-1}a_{n,k}\left\vert S_{\frac{%
\alpha k}{2}}f\left( x\right) -f\left( x\right) \right\vert ^{q}\right\}
^{1/q}+\left\{ \sum_{m=\left[ c\right] }^{\infty
}\sum_{k=2^{m}}^{2^{m+1}-1}a_{n,k}\left\vert S_{\frac{\alpha k}{2}}f\left(
x\right) -f\left( x\right) \right\vert ^{q}\right\} ^{1/q} \\
&=&I_{1}\left( x\right) +I_{2}\left( x\right) .
\end{eqnarray*}%
Using Proposition 4 and denoting the left hand side of the inequality from
its by $F_{n}$ ,i.e. $F_{n}=w_{x}\left( \frac{\pi }{n+1}\right) +E_{\alpha
n/2}\left( f\right) _{S^{p}},$ we get 
\begin{eqnarray*}
I_{1}\left( x\right) &\leq &\left\{ \sum_{k=0}^{2^{\left[ c\right]
}-1}a_{n,k}\frac{k/2+1}{k/2+1}\sum\limits_{l=k/2}^{k}\left\vert S_{\frac{%
\alpha l}{2}}f\left( x\right) -f\left( x\right) \right\vert ^{q}\right\}
^{1/q} \\
&\leq &\left\{ 2^{\left[ c\right] }\sum_{k=0}^{2^{\left[ c\right] }-1}a_{n,k}%
\frac{1}{k/2+1}\sum\limits_{l=k/2}^{k}\left\vert S_{\frac{\alpha l}{2}%
}f\left( x\right) -f\left( x\right) \right\vert ^{q}\right\} ^{1/q} \\
&\ll &\left\{ \sum_{k=0}^{2^{\left[ c\right] }-1}a_{n,k}F_{k/2}^{q}\right\}
^{1/q}.
\end{eqnarray*}%
By partial summation, our Proposition 4 gives 
\begin{eqnarray*}
I_{2}^{q}\left( x\right) &=&\sum_{m=[c]}^{\infty }\left[
\sum_{k=2^{m}}^{2^{m+1}-2}\left( a_{n,k}-a_{n,k+1}\right)
\sum_{l=2^{m}}^{k}\left\vert S_{\frac{\alpha l}{2}}f\left( x\right) -f\left(
x\right) \right\vert ^{q}\right. \\
&&\left. +a_{n,2^{m+1}-1}\sum_{l=2^{m}}^{2^{m+1}-1}\left\vert S_{\frac{%
\alpha l}{2}}f\left( x\right) -f\left( x\right) \right\vert ^{q}\right]
\end{eqnarray*}%
\begin{eqnarray*}
&\ll &\sum_{m=[c]}^{\infty }\left[ 2^{m}\sum_{k=2^{m}}^{2^{m+1}-2}\left\vert
a_{n,k}-a_{n,k+1}\right\vert F_{\alpha 2^{m}/2}^{q}\right. \\
&&\left. +2^{m}a_{n,2^{m+1}-1}F_{\alpha 2^{m}/2}^{q}\right] \\
&=&\sum_{m=[c]}^{\infty }2^{m}F_{\alpha 2^{m}/2}^{q}\left[
\sum_{k=2^{m}}^{2^{m+1}-2}\left\vert a_{n,k}-a_{n,k+1}\right\vert
+a_{n,2^{m+1}-1}\right] .
\end{eqnarray*}%
Since (\ref{5}) holds, we have 
\begin{eqnarray*}
&&a_{n,s+1}-a_{n,r} \\
&\leq &\left\vert a_{n,r}-a_{n,s+1}\right\vert \leq \sum_{k=r}^{s}\left\vert
a_{n,k}-a_{n,k+1}\right\vert \\
&\leq &\sum_{k=2^{m}}^{2^{m+1}-2}\left\vert a_{n,k}-a_{n,k+1}\right\vert \ll
\sum\limits_{k=[2^{m}/c]}^{[c2^{m}]}\frac{a_{n,k}}{k}\text{ \ \ }\left(
2\leq 2^{m}\leq r\leq s\leq 2^{m+1}-2\right) ,
\end{eqnarray*}%
whence 
\begin{equation*}
a_{n,s+1}\ll a_{n,r}+\sum\limits_{k=[2^{m}/c]}^{[c2^{m}]}\frac{a_{n,k}}{k}%
\text{ \ }\left( 2\leq 2^{m}\leq r\leq s\leq 2^{m+1}-2\right)
\end{equation*}%
and 
\begin{eqnarray*}
2^{m}a_{n,2^{m+1}-1} &=&\frac{2^{m}}{2^{m}-1}%
\sum_{r=2^{m}}^{2^{m+1}-2}a_{n,2^{m+1}-1} \\
&\ll &\sum_{r=2^{m}}^{2^{m+1}-2}\left(
a_{n,r}+\sum\limits_{k=[2^{m}/c]}^{[c2^{m}]}\frac{a_{n,k}}{k}\right) \\
&\ll
&\sum_{r=2^{m}}^{2^{m+1}-1}a_{n,r}+2^{m}\sum\limits_{k=[2^{m}/c]}^{[c2^{m}]}%
\frac{a_{n,k}}{k}.
\end{eqnarray*}%
Thus 
\begin{equation*}
I_{2}^{q}\left( x\right) \ll \sum_{m=[c]}^{\infty }\left\{ 2^{m}F_{\alpha
2^{m}/2}^{q}\sum\limits_{k=[2^{m}/c]}^{[c2^{m}]}\frac{a_{n,k}}{k}+F_{\alpha
2^{m}/2}^{q}\sum_{k=2^{m}}^{2^{m+1}-1}a_{n,k}\right\} .
\end{equation*}

Finally, by elementary calculations we get 
\begin{eqnarray*}
I_{2}^{q}\left( x\right) &\ll &\sum_{m=[c]}^{\infty }\left\{ 2^{m}F_{\alpha
2^{m}/2}^{q}\sum\limits_{k=2^{m-\left[ c\right] }}^{2^{m+\left[ c\right] }}%
\frac{a_{n,k}}{k}+F_{\alpha
2^{m}/2}^{q}\sum_{k=2^{m}}^{2^{m+1}}a_{n,k}\right\} \\
&\ll &\sum_{m=[c]}^{\infty }F_{\alpha 2^{m}/2}^{q}\sum\limits_{k=2^{m-\left[
c\right] }}^{2^{m+\left[ c\right] }}a_{n,k} \\
&=&\sum_{m=[c]}^{\infty }F_{\alpha 2^{m}/2}^{q}\sum\limits_{k=2^{m-\left[ c%
\right] }}^{2^{m}-1}a_{n,k}+\sum_{m=[c]}^{\infty }F_{\alpha
2^{m}/2}^{q}\sum\limits_{k=2^{m}}^{2^{m+\left[ c\right] }}a_{n,k}
\end{eqnarray*}%
\begin{equation*}
\ll \sum_{m=[c]}^{\infty }\sum\limits_{k=2^{m-\left[ c\right]
}}^{2^{m}-1}a_{n,k}F_{\alpha k/2}^{q}+\sum_{m=[c]}^{\infty
}\sum\limits_{k=2^{m}}^{2^{m+\left[ c\right] }}a_{n,k}F_{\frac{\alpha k}{%
2^{1+\left[ c\right] }}}^{q}
\end{equation*}%
\begin{equation*}
=\sum_{m=[c]}^{\infty }\sum\limits_{k=2^{m-\left[ c\right]
}}^{2^{m}-1}a_{n,k}F_{\alpha k/2}^{q}+\sum_{m=[c]}^{\infty
}\sum\limits_{k=2^{m}}^{2^{m+\left[ c\right] }-1}a_{n,k}F_{\frac{\alpha k}{%
2^{1+\left[ c\right] }}}^{q}+\sum_{m=[c]}^{\infty }F_{\frac{\alpha 2^{m}}{2}%
}^{q}a_{n,2^{m+\left[ c\right] }}
\end{equation*}%
\begin{eqnarray*}
&=&\sum_{m=[c]}^{\infty }\sum\limits_{r=1}^{\left[ c\right]
}\sum\limits_{k=2^{m-r}}^{2^{m-r+1}-1}a_{n,k}F_{\alpha
k/2}^{q}+\sum_{m=[c]}^{\infty }\sum\limits_{r=0}^{\left[ c\right]
-1}\sum\limits_{k=2^{m+r}}^{2^{m+r+1}-1}a_{n,k}F_{\frac{\alpha k}{2^{1+%
\left[ c\right] }}}^{q} \\
&&+\sum_{m=[c]}^{\infty }F_{\frac{\alpha 2^{m}}{2}}^{q}a_{n,2^{m+\left[ c%
\right] }}
\end{eqnarray*}%
\begin{eqnarray*}
&\leq &\sum\limits_{r=1}^{\left[ c\right] }\sum\limits_{k=2^{\left[ c%
\right] -r}}^{\infty }a_{n,k}F_{\alpha k/2}^{q}+\sum\limits_{r=0}^{\left[ c%
\right] -1}\sum\limits_{k=2^{\left[ c\right] +r}}^{\infty }a_{n,k}F_{\frac{%
\alpha k}{2^{1+\left[ c\right] }}}^{q}+\sum\limits_{k=2^{2\left[ c\right]
}}^{\infty }a_{n,k}F_{\frac{\alpha k}{2^{1+\left[ c\right] }}}^{q} \\
&\ll &\sum\limits_{k=0}^{\infty }a_{n,k}F_{\frac{\alpha k}{2^{1+\left[ c%
\right] }}}^{q}.
\end{eqnarray*}

Thus we obtain the desired result. $\square $

\subsection{Proof of Theorem 6}

If $\left( a_{n,k}\right) _{k=0}^{\infty }\in MS$ then $\left(
a_{n,k}\right) _{k=0}^{\infty }\in GM\left( _{2}\beta \right) $ and using
Theorem 5 we obtain 
\begin{eqnarray*}
H_{n,A,\gamma }^{q}f\left( x\right) &\leq &\left\{ \sum_{k=0}^{\infty
}a_{n,k}\left[ w_{x}(\frac{\pi }{k+1})\right] ^{q}\right\} ^{1/q} \\
&&+\left\{ \sum_{k=0}^{\infty }\sum\limits_{m=k2^{\left[ c\right]
}}^{\left( k+1\right) 2^{\left[ c\right] }-1}a_{n,m}\left[ E_{\frac{\alpha m%
}{2^{1+\left[ c\right] }}}\left( f\right) _{S^{p}}\right] ^{q}\right\} ^{1/q}
\end{eqnarray*}%
\begin{eqnarray*}
&\leq &\left\{ \sum_{k=0}^{\infty }a_{n,k}\left[ w_{x}(\frac{\pi }{k+1})%
\right] ^{q}\right\} ^{1/q} \\
&&+\left\{ \sum_{k=0}^{\infty }\sum\limits_{m=k2^{\left[ c\right]
}}^{\left( k+1\right) 2^{\left[ c\right] }-1}a_{n,m}\left[ E_{\frac{\alpha k%
}{2}}\left( f\right) _{S^{p}}\right] ^{q}\right\} ^{1/q}
\end{eqnarray*}%
\begin{eqnarray*}
&\leq &\left\{ \sum_{k=0}^{\infty }a_{n,k}\left[ w_{x}(\frac{\pi }{k+1})%
\right] ^{q}\right\} ^{1/q} \\
&&+\left\{ \sum_{k=0}^{\infty }2^{\left[ c\right] }a_{n,k2^{\left[ c\right]
}}\left[ E_{\frac{\alpha k}{2}}\left( f\right) _{S^{p}}\right] ^{q}\right\}
^{1/q} \\
&\ll &\left\{ \sum_{k=0}^{\infty }a_{n,k}\left[ w_{x}(\frac{\pi }{k+1})+E_{%
\frac{\alpha k}{2}}\left( f\right) _{S^{p}}\right] ^{q}\right\} ^{1/q}
\end{eqnarray*}%
This ends our proof. $\square $

\end{document}